\documentclass[preprint,12pt]{elsarticle}
\usepackage[symbol]{footmisc}

\usepackage{amsmath}
\usepackage{amsfonts}
\usepackage{amssymb}
\usepackage{amsthm}
\usepackage{hyperref}
\usepackage{url}
\usepackage{indentfirst}
\usepackage{tikz}
\usepackage{booktabs}
\usepackage{pgfplots}
\pgfplotsset{width=\textwidth, compat=newest}
\usetikzlibrary{shapes.geometric}
\usepackage{filecontents}
\usepackage{booktabs}
\usepackage[noend]{algpseudocode}
\usepackage[margin=1.5in]{geometry}    
\usepackage{graphicx}
\usetikzlibrary{shapes,positioning}
\usepackage[utf8]{inputenc}
\usepackage{hyperref}
\usepackage{caption}
\usepackage{epigraph}
\usepackage{diagbox}

\usepackage[justification=centering]{caption} 
\usepackage{graphicx}
\usepackage[export]{adjustbox}
\usepackage{caption, subcaption}

\usepackage{chngcntr}
\counterwithout{remark}{section}
\counterwithout{figure}{section}
\counterwithout{table}{section}
\counterwithout{theorem}{section}
\counterwithout{lemma}{section}
\counterwithout{proposition}{section}
\journal{the Notices of AMS}

\begin{document}
\begin{frontmatter}

\title{Supercentenarian paradox}

\author{Vladimir Gurvich}
\ead{vgurvich@hse.ru and vladimir.gurvich@gmail.com}
\address{National Research University Higher School of Economics, Moscow, Russia}
\ead{vladimir.gurvich@gmail.com}
\address{RUTCOR, Rutgers University, Piscataway, New Jersey, United States}

\author{Mariya Naumova}
\ead{mnaumova@math.rutgers.edu}
\address{Rutgers Business School, Rutgers University, Piscataway, New Jersey, United States}

\end{frontmatter}

\epigraph{And the LORD said, My spirit shall not always strive with man,
for that he also is flesh: yet his days shall be one hundred and twenty years.}{Genesis 6:3}
\epigraph{What should you wish for the 120th birthday? Have a nice day!}{A corollary}

\epigraph{Yes, man is mortal, but that would be only half the trouble. The worst of it is that he's sometimes unexpectedly mortal—there's the trick!}{Mikhail Bulgakov, The Master and Margarita}

\section{To bet or not to bet: that is the question }
\label{s1}

Consider the following statement:

\centerline{
$B(t, \Delta t)$:  a  $t$-year-old person NN will survive another $\Delta t$  years,}
\noindent
where $t, \Delta t \in \mathbb{R_+}$  are nonnegative real numbers.

We only know that $NN$  is  $t$  years old and nothing else, in particular, gender, origin, health conditions, etc. are not known.
Should we bet that  $B(t, \Delta t)$  will hold?
It seems that our odds are very good for any  $t$  provided
$\Delta t$ is small enough, say, $1 / 365$ (that is, one day).
However, this is not that obvious, and it depends on the life-time probability distribution.

Let  $F(t)$ be the cumulative distribution function of the lifetime of a person, where $t$ is in years, and define $\Phi(t) = 1 - F(t)$.

The following cases are possible:
\begin{enumerate}
    \item
    $\Phi(t)=0$ whenever $t \ge t_0$ for some 
    $t_0 \in \mathbb{R_+}$.

    In this scenario, if $t \in [t_0 - \Delta t, t_0]$, the bet is lost with probability 1. 
    
    \item
    $\Phi(t) \rightarrow 0$ faster than an exponential distribution.
    
    Then, for each fixed $\Delta t > 0$, the bet is lost with probability as close to 1 as desired if $t$ is large enough.
    \item
   $\Phi(t) \rightarrow 0$ as fast as an exponential distribution.

   Then for large $t$ the chances to win the bet depend only on  $\Delta t$ and not on $t$.
   \item
   $\Phi(t) \rightarrow 0$ slower than an exponential distribution.
   
   In this scenario, in agreement with intuition, for any fixed $t$, the chances to win tend to $1$ as $\Delta t \rightarrow 0$. 
\end{enumerate}

The last two cases are not realistic as statistical literature suggests; see below.

For example, $normal(\mu, \sigma^2)$ distribution corresponds to case 2, as

\begin{equation}
\begin{aligned}
\displaystyle \lim_{t \rightarrow \infty} \frac{\Phi(t+\Delta t)}{\Phi(t)} &= \displaystyle \lim_{t \rightarrow \infty} \frac{1-\displaystyle \int_{- \infty}^{t+\Delta t}\frac{1}{\sqrt{2\pi}\sigma}e^{-\frac{(x-\mu)^2}{2\sigma^2}}dx}{1-\displaystyle \int_{- \infty}^{t}\frac{1}{\sqrt{2\pi}\sigma}e^{-\frac{(x-\mu)^2}{2\sigma^2}}dx}\\
 & = \displaystyle \lim_{t \rightarrow \infty} \frac{e^{-\frac{(t+\Delta t-\mu)^2}{2\sigma^2}}}{e^{-\frac{(t-\mu)^2}{2\sigma^2}}}\\
 & = \displaystyle \lim_{t \rightarrow \infty} e^{-\frac{2t\Delta t +(\Delta t)^2-2\mu\Delta t}{2\sigma^2}}\\
 & = 0 \hspace{0.2in} \text{for any fixed } \Delta t>0.\nonumber
\end{aligned}
\end{equation}

Clearly, all distributions with negative excess kurtosis (platykurtic distributions), that is, the distributions with thinner tails compared to the normal distribution, satisfy this property.

In demography and actuarial sciences, the Gompertz distribution \cite{Gom1825} is often applied to describe the distribution of adult lifespans \cite{PHG2001}, \cite{Spi2016}, \cite{GK2021}, \cite{KGB2021}, \cite{Kul2021}; see Section 2 for more details.

The probability density function of the Gompertz distribution is given by

\begin{equation}
\displaystyle f\left(x;\eta ,b\right)=b\eta \exp \left(\eta +bx-\eta e^{bx}\right){\text{  for }}x\geq 0,\,\nonumber
\end{equation}
where $b>0$ is the scale parameter and $\eta >0$ is the shape parameter.

The cumulative distribution function is
\begin{equation}
F\left(x;\eta ,b\right)=1-\exp \left(-\eta \left(e^{bx}-1\right)\right).\nonumber
\end{equation}

We find that 
\begin{equation}
p = \frac{\Phi(t+\Delta t)}{\Phi(t)} = \frac{e^{-\eta(e^{b(t+\Delta t)}-1)}}{e^{-\eta(e^{bt}-1)}}=e^{-\eta e^{bt}(e^{b\Delta t}-1)}.\nonumber
\end{equation}

The minimal time $t$ such that the probability to win the bet is at least $p$ is 
$$t = \frac{1}{b}\ln\left(\frac{-\ln p}{\eta (e^{b\Delta t}-1)}\right),$$
and the combination of parameters $\eta$ and $b$ that yields the distribution mode of 75 \footnote[1]{For this example, we choose the mode of the distribution to be 75 years. Note that the life expectancy at birth in the US depending on state for 2019 was in statistical range of 74.4-80.9 years; see \url{https://www.cdc.gov/nchs/data-visualization/state-life-expectancy/index_2019.htm.}}
is $b=-\frac{1}{75}\ln \eta$, $\eta \in (0, 1)$. Statisticians fit the parameters $b$ and $\eta$ of the distribution to a series of data based on various factors (population average age, geographic location, etc). The values of $t$ for various $\Delta t$ with 50$\%$ chance of winning the bet are shown in Table \ref{t1}. 


\begin{table}
\centering
\caption{The values of $t$ for various $\Delta t$ with 50$\%$ chance of winning the bet in case of several Gompertz distributions with mode $75$}
\small
\begin{tabular}{|l|cccccc|}
\hline
 \diagbox[width=\dimexpr \textwidth/8+2\tabcolsep\relax, height=1cm]{$\Delta t$}{$\eta$}&$e^{-75/3}$&$e^{-75/4}$&$e^{-75/5}$&$e^{-75/6}$&$e^{-75/7}$&$e^{-75/8}$\\
\hline
\hline
1 year &76.68&	78.57&	80.71&	83.04&	85.55&	88.20\\
\hline
1 month &84.61&	88.98&	93.60&	98.42	&103.41&	108.54\\
\hline
1 day&94.89&	102.68&	110.71&	118.95&	127.35	&135.90\\
\hline
1 hour&104.43&	115.39&	126.60&	138.02&	149.60&	161.33\\
\hline
1 minute&116.71&	131.77&	147.08&	162.59&	178.26&	194.08\\
\hline
1 second& 129.00&	148.15&	167.55&	187.15&	206.92&	226.84\\
\hline
\hline
\end{tabular}
\label{t1}
\end{table}

\begin{figure}[h]
\begin{tikzpicture}
\begin{axis}[
    legend style={
    legend pos=north east,
    legend style={row sep=1.5pt}
},
    width=14cm,
    axis x line = bottom,
    axis y line = left,
    ylabel={t},
    xlabel={\(\Delta t\)},
    xtick={1/12, 1/365},
    xticklabels={$\frac{1}{12}$, $\frac{1}{365}$},
    every axis x label/.style={at={(current axis.right of origin)},anchor=west},
        every axis y label/.style={at={(current axis.north west)}, above=0.5mm}
    ]
    \addplot [red, dotted, 
        domain=0:1/10,
        samples=1001,
        very thick
        ]
        {3*(ln(ln(2))+75/3-ln(exp(x/3)-1))}; 
    \addlegendentry{$\eta = e^{-75/3}$}
    \addplot [blue, dotted,
        domain=0:1/10,
        samples=1001,
        very thick
        ]
        {4*(ln(ln(2))+75/4-ln(exp(x/4)-1))};
    \addlegendentry{$\eta = e^{-75/4}$}
    \addplot [brown, dotted,
        domain=0:1/10,
        samples=1001,
        very thick
        ]
        {5*(ln(ln(2))+75/5-ln(exp(x/5)-1))};
    \addlegendentry{$\eta = e^{-75/5}$}
    \addplot [black, dotted,
        domain=0:1/10,
        samples=1001,
        very thick
        ]
        {6*(ln(ln(2))+75/6-ln(exp(x/6)-1))};
    \addlegendentry{$\eta = e^{-75/6}$}
    \addplot [magenta, dotted,
        domain=0:1/10,
        samples=1001,
        very thick
        ]
        {7*(ln(ln(2))+75/7-ln(exp(x/7)-1))};
    \addlegendentry{$\eta = e^{-75/7}$}
    \addplot [green, dotted,
        domain=0:1/10,
        samples=1001,
        very thick
        ]
        {8*(ln(ln(2))+75/8-ln(exp(x/8)-1))};
    \addlegendentry{$\eta = e^{-75/8}$}
    \end{axis}
\end{tikzpicture}
\caption{Relation between $\Delta t$ and $t$, for $\Delta t < 1/12$ (1 month), with 50$\%$ chance of winning the bet in case of several Gompertz distributions with mode $75$}
\label{f2}
\end{figure}
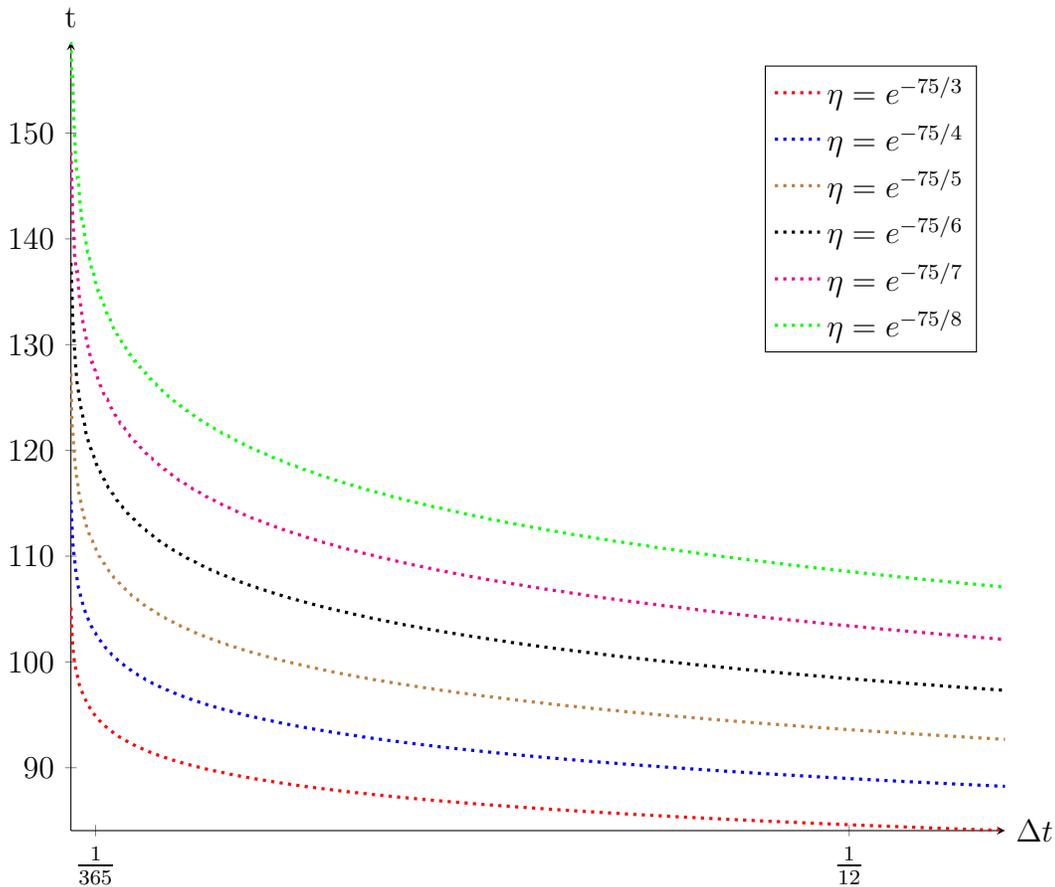

Figure \ref{f2} demonstrates how high $t$ values are for small values of $\Delta t$ (one month or less). 
For example, if $b=1/3$, $\eta = e^{-75/3}$, then $\Delta t = 1/(365\cdot 24 \cdot 60 \cdot 60)$ (one second) requires $t \approx 129$ (years), but on Earth there is no person of such age \cite{DMV2016}.

\section{Actuary mortality models}
Life actuaries in the US use the mortality tables published by the Society of Actuaries (SOA; see Mort.soa.org), which are based on the collected data. The Continuous Mortality Investigation (CMI) in the UK provides projections for the expected changes to those tables. Most insurance companies use these tables, adjusting them to their needs, while others use their own historical data to cover probabilities down to a personal level, including variables such as location, income, and health conditions; see \url{https://www.actuaries.org.uk/learn-and-develop/continuous-mortality-investigation}. 

Meanwhile, life contingencies textbooks often have hypothetical life tables based on Gompertz and Gompertz-Makeham models \cite{PHG2001}, \cite{MLV2012}. Despite their age, the models, also known as the law of mortality, are still widely used in the fields of demographics, biology, actuarial science and even financial engineering \cite{Mil2020}, \cite{Kul2021}. Modern actuaries incorporate the Gompertz distribution to represent important actuarial expressions in closed form via the Gamma function, which makes representation especially convenient and popular in the annuity literature \cite{Mil2020}.

The more advanced Gompertz–Makeham model assumes that the human death rate is the sum of two components: an age-dependent term (the Gompertz function \cite{Gom1825}), which increases exponentially with age, and an age-independent term (the Makeham function \cite{Mak1860}). For low mortality countries the latter term is often negligible and the formula simplifies to the Gompertz model discussed in Section 1. While, based on the world data, the Makeham term has been in decline since the 1950s, the Gompertz component is surprisingly stable.

\medskip

One can verify that for the above-mentioned distributions, $Pr(B(t, \Delta t)) \rightarrow 0$  as  $t \rightarrow \infty$, for any fixed  $\Delta t > 0$. Hence, for arbitrarily small positive  $\Delta t $ and $\epsilon$ there exists a sufficiently large $t$  such that  $Pr(B(t, \Delta t))  < \epsilon$, which means that we should not bet...  in theory. However, in practice we can bet safely, because for any fixed positive $\Delta t$ and $\epsilon$, a very large $t$ is required to satisfy the inequality  $Pr(B(t, \Delta t))  < \epsilon$.

For example, $\Delta t = 1/365$ and $\epsilon = 1/2$  may require  $t > 125$  years for some  typical distributions  $F$  considered in the literature. Yet, on Earth there is no person of such age. Thus, our odds are good because the chosen testee  NN is not old enough, or the bet cannot be made for technical (or, more precisely, statistical) reasons - absence of a testee.

However, there is not much statistical data about supercentenarians due to their scarceness \footnote{According to Guinness World Records, the oldest person currently living is Kane Tanaka (Japan, b. 2 January 1903) aged 119 years and 18 days, as verified on 20 January 2022; the only person verified to have lived beyond the age of 120 is Jeanne Louise Calment (France, 1875-1997).}. Therefore, the tail distribution of the human lifetime can be fitted based on an extrapolation only. Yet, it is well-known that extrapolation is not reliable. As Mark Twain noticed:

\medskip

\textit{In the space of one hundred and seventy-six years the Lower Mississippi has shortened itself two hundred and forty-two miles. That is an average of a trifle over one mile and a third per year. Therefore, any calm person, who is not blind or idiotic, can see that in the Old Oolitic Silurian Period, just a million years ago next November, the Lower Mississippi River was upwards of one million three hundred thousand miles long, and stuck out over the Gulf of Mexico like a fishing-rod. And by the same token any person can see that seven hundred and forty-two years from now the Lower Mississippi will be only a mile and three-quarters long, and Cairo and New Orleans will have joined their streets together, and be plodding comfortably along under a single mayor and a mutual board of aldermen. There is something fascinating about science. One gets such wholesale returns of conjecture out of such a trifling investment of fact.} (Mark Twain, Life in the Mississippi).

The Supercentenarian Paradox is similar to the famous St.Petersburg Paradox, which we briefly review in the following section.

\section{Analogy with St. Petersburg's Paradox}

The St. Petersburg paradox was named after one of the leading scientific journals of the eighteenth century Commentarii Academiae Scientiarum Imperialis Petropolitanae (Papers of the Imperial Academy of Sciences in Petersburg), in which Daniel Bernoulli published a paper entitled “Specimen Theoriae Novae de Mensura Sortis” (“Exposition of a New Theory on the Measurement of Risk”) in 1738 \cite{Ber1738}. The problem was introduced by Nicolas Bernoulli \cite{Plo1993} who stated it in a letter to Pierre Raymond de Montmort on September 9, 1713  \cite{Mon1713}, \cite{Eve1990}.  

 The paradox appears from the game that is played as follows:

\textit{A fair coin is flipped until it comes up heads the first time.
Let  $\; t = 0,1,2, \ldots$ denote the number of preceding tails.
Then,  the gambler  wins  $2^{t+1}$  dollars. 
What would be a fair compensation to pay the casino for entering the game?}

Assuming that the casino has unlimited resources,  the expected gain of the gambler is

$${\displaystyle {\begin{aligned}Gain&= \displaystyle \sum_{t=0}^\infty 2^{t+1} (1/2)^{t+1} = 
{\frac {1}{2}}\cdot 2+{\frac {1}{4}}\cdot 4+{\frac {1}{8}}\cdot 8+{\frac {1}{16}}\cdot 16+\cdots \\&=1+1+1+1+\cdots \\&=\infty \,.\end{aligned}}}$$

Thus, in theory it is rational for the gambler to pay the casino
an arbitrary large amount for a single opportunity to play the St. Petersburg game. Yet, almost certainly the gambler will win a modest amount. Even $\$20$ compensation looks too generous.
Nicolas Bernoulli conjectured that people will neglect unlikely events \cite{Mon1713}.

\medskip

Similarly, in the Supercentenarian paradox, statement
$B(t, \Delta t)$  holds for arbitrarily small  $\Delta t$ with arbitrary small probability, 
provided the age $t$  is sufficiently large.
In both cases the contribution of the very large values of $t$
(that is, very long sequences of tails and very old people) is decisive, but only in theory. 
In practice such values of  $t$  are too unlikely.

Interestingly, already in 1777, Buffon  \cite{Buf1777} 
mentioned certain connections between the St. Petersburg Paradox  and  $B(t, \Delta t)$.  As cited in  \cite{Dut1988},  
\medskip

\textit{[Buffon] notes that a fifty-six year old man, believing his health to be good, would disregard the probability that he would die within twenty-four hours, although mortality tables indicate that the odds against his dying in this period are only 10189 to 1. Buffon thus takes a probability of 1/10,000 or less for an event as a probability which may be disregarded.}

\end{document}